\documentclass[12pt]{amsart}
\usepackage{cases}
\usepackage{txfonts}
\usepackage{latexsym}
\usepackage{amsthm}
\usepackage{mathrsfs}
\usepackage{amssymb, amsmath}
\usepackage{dsfont}
\def\opn#1#2{\def#1{\operatorname{#2}}} 
\opn\chara{char} \opn\length{\ell}
\opn\projdim{proj\,dim} \opn\injdim{inj\,dim} \opn\rank{rank}
\opn\depth{depth} \opn\grade{grade} \opn\height{height}
\opn\embdim{emb\,dim} \opn\codim{codim}

\opn\Tr{Tr} \opn\bigrank{big\,rank}
\opn\superheight{superheight}\opn\lcm{lcm}
\opn\trdeg{tr\,deg}%
\opn\reg{reg} \opn\lreg{lreg}
\opn\div{div} \opn\Div{Div} \opn\cl{cl} \opn\Cl{Cl}
\opn\Spec{Spec} \opn\Supp{Supp} \opn\supp{supp} \opn\Sing{Sing}
\opn\Ass{Ass}
\opn\Ann{Ann} \opn\Rad{Rad} \opn\Soc{Soc}
\opn\Ker{Ker} \opn\Coker{Coker} \opn\Im{Im} \opn\Hom{Hom}
\opn\Tor{Tor} \opn\Ext{Ext} \opn\End{End} \opn\Aut{Aut} \opn\id{id}

\opn\nat{nat}
\opn\pff{pf}
\opn\Pf{Pf} \opn\GL{GL} \opn\SL{SL} \opn\mod{mod} \opn\ord{ord}
\opn\aff{aff} \opn\con{conv} \opn\relint{relint} \opn\st{st}
\opn\lk{lk} \opn\cn{cn} \opn\core{core} \opn\vol{vol}
\opn\gr{gr}

\def\pot#1#2{#1[\kern-0.28ex[#2]\kern-0.28ex]}

\opn\dirlim{\underrightarrow{\lim}}
\opn\invlim{\underleftarrow{\lim}}
\newtheorem{theorem}{Theorem}[section]
\newtheorem{corollary}[theorem]{Corollary}
\newtheorem{lemma}[theorem]{Lemma}
\newtheorem{conjecture}[theorem]{Conjecture}
\newtheorem{definition}[theorem]{Definition}
\newtheorem{question}[theorem]{Question}
\def\qed{\ifhmode\textqed\fi
   \ifmmode\ifinner\quad\qedsymbol\else\dispqed\fi\fi}
\def\textqed{\unskip\nobreak\penalty50
    \hskip2em\hbox{}\nobreak\hfil\qedsymbol
    \parfillskip=0pt \finalhyphendemerits=0}
\def\dispqed{\rlap{\qquad\qedsymbol}}

\opn\ini{in} \opn\inm{inm} \opn\Sym{Sym}

\begin{document}
\title{On the set of the difference of primes}
\author{Wen Huang}
\address {School of Mathematical Sciences\\ University of Science and Technology of China\\ Hefei 230026, P. R.
China.}
\email {wenh@mail.ustc.edu.cn}

\author{XiaoSheng Wu}
\address {School of Mathematics, Heifei University of Technology, Hefei 230009,
P. R. China}
\email {xswu@amss.ac.cn}
\date{}


\begin{abstract} In this work we prove that the set of the difference of primes is a $\Delta_r^*$-set. The work is based on the recent dramatic new developments in the study of bounded gaps between primes, reached by Zhang, Maynard and Tao.

\end{abstract}
\maketitle

\section{Introduction}
\label{sec1}
In the present work we investigate on how ``large'' the set of the difference of primes. In 1905, Maillet conjectured in \cite{mai} that the set of the difference of primes should have the ``largest'' form, which contains all even numbers.
\begin{conjecture}
\label{conjecture1}
Every even number is the difference of two primes.
\end{conjecture}
Actually, before Maillet's conjecture, there were two stronger forms of the conjecture. One was due to Kronecker \cite{kro} in 1901:
\begin{conjecture}
\label{conjecture2}
Every even number can be expressed in infinitely many ways as the difference of two primes.
\end{conjecture}
Another is formulated by Polignac \cite{pol} in 1849, which has the most general form as
\begin{conjecture}
\label{conjecture3}
Every even number can be written in infinitely many ways as the difference of two consecutive primes.
\end{conjecture}
It is easy to see that the twin prime conjecture is about the lower bound for the set of the difference of primes in Conjecture \ref{conjecture2} or \ref{conjecture3}. Recently, based on the GPY\cite{gol} sieve method, Zhang \cite{zha} made a breakthrough and proved that there exists an even number not more than $7\times10^7$ which can be expressed in infinitely many ways as the difference of two primes. Soon after, Maynard \cite{may} and Tao reduced the limit of such even number to not more than 600. The best known result now is not more than 246, to see \cite{pol}. Assume that the primes have level of distribution $\theta$ for every $\theta<1$, then the best known result now is 12 proved by Maynard \cite{may}. Here, for some given $\theta>0$, we say the primes have `level of distribution $\theta$' if, for any $W>0$, we have
\begin{align}
   \sum_{q\le x^\theta}\max_{(a,q)=1}\Big|\pi(x;q,a)-\frac{\pi(x)}{\phi(q)}\Big|\ll_W\frac{x}{(\log x)^W}.
\end{align}

Another important aspect of these conjectures is how ``large'' the set of the difference of primes. In combinatorial number theory, as well as in dynamics, there are various notions of ``large'' sets of integers. Some familiar notions are those of sets of positive (upper) density, syndetic sets, thick set, return-time sets, sets of recurrence, Bohr sets, Nil$_d$ Bohr$_0$-set, piecewise-Bohr sets and strongly piecewise-Bohr sets. We will give some basic definitions and elementary considerations of these notions in section \ref{sec2}.

Let $D$ denote the set of even numbers that can be expressed in infinitely many ways as the difference of two primes. Based on recent breakthrough on the twin prime conjecture, Pintz proved in \cite{pin} that
\begin{theorem}
 There exists an ineffective constant $C'$ such that every interval of type $[M, M+C']$ contains at least one even numbers that can be expressed in infinitely many ways as the difference of two primes, that is to say $D$ is a syndetic set.
\end{theorem}
Somewhat later, using a
different method, Granville, Kane, Koukoulopoulos and Lemke Oliver obtained the same result in \cite{gra}.

\begin{definition}
If $S$ is a non-empty subset of $\mathbb{N}$, define the difference set $\Delta(S)$ by
\begin{align}
   \Delta(S)=(S-S)\cap\mathbb{N}=\{b-a:a\in S, b\in S, b> a\}.
\end{align}
If $A$ is a subset of $\mathbb{N}$, $A$ is a $\Delta_{r}^*$-set if $A\cap\Delta(S)\neq\emptyset$ for every subset $S$ of $\mathbb{N}$ with $|S|=r$; $A$ is a $\Delta^*$-set if $A\cap\Delta(S)\neq\emptyset$ for every infinite subset $S$ of $\mathbb{N}$.
\end{definition}

In this paper, we prove that $D$ is a ``larger'' set than a syndetic set.
\begin{theorem}
\label{theoremhw}
Let $D$=\{$d$: $d$ can be expressed in infinitely many ways as the difference of two primes \}, then $D$ is a $\Delta_{r}^*$-set for $r\ge721$. If assume that the primes have level of distribution $\theta$ for every $\theta<1$, we have $D$ is a $\Delta_{r}^*$-set for $r\ge19$.
\end{theorem}
Actually, we obtain the following inequality for the lower bound of $r$ in Theorem \ref{theoremhw}.
\begin{align}
     r\ge C\prod_{p\le C}\big(1-\frac1p\big)^{-1},
\end{align}
where $C$ is the lower bound of the length of admissible $k$-tuple of integers in Zhang-Maynard-Tao's theorem, to see section \ref{sec3} in the following.

It is proved by Green and Tao in \cite{gre} that the primes contain arbitrarily long arithmetic progressions, but we know little about the distribution of the common difference for these arithmetic progressions. Let $\mathbb{P}$ denote the set of all primes and
\begin{align}
  C_d(\mathbb{P})=\{n\in\mathbb{N}:\mathbb{P}\cap(\mathbb{P}-n) \cap(\mathbb{P}-2n)\cap\cdots\cap(\mathbb{P}-dn)\neq\emptyset\}.
\end{align}
Huang, Shao and Ye asked the following question in their work \cite{hua}
\begin{question}
\label{que}
Is $C_d(\mathbb{P})$ a Nil$_d$ Bohr$_0$-set~?
\end{question}

When take $d=1$, we can see that Question \ref{que} actually asks whether $D$ is a Bohr$_0$-set. From Theorem \ref{theoremhw}, we have the following corollary for this question
\begin{corollary}
\label{cor}
Let $D$=\{$d$: $d$ can be expressed in infinitely many ways as the difference of two primes \}, then $D$ is a strongly piecewise-$\text{Bohr}_0$ set.
\end{corollary}

\section{some notations}
\label{sec2}

We begin with basic definitions and elementary considerations for some notions. A set $S\subset \mathbb{N}$ is called \emph{syndetic} if for some finite subset $F\subset\mathbb{N}$
\begin{align}
\bigcup_{n\in F}\Big(S-n\Big)=\mathbb{N},
\end{align}
where $S-n=\{m\in\mathbb{N}: m+n\in S\}$. In other words, $S$ is a syndetic set if it has bounded gaps, which means that there is an integer $k$ such that $\{a, a+1, a+2, \cdots, a+k\}\cap S\neq\emptyset$ for any $a\in\mathbb{N}$.

A set $A\subset \mathbb{N}$ is called a \emph{thick} set if it contains arbitrarily long intervals. That is, for every $m\in\mathbb{N}$, there is some $n\in\mathbb{N}$ such that $\{n, n+1, n+2, \cdots, n+m\}\subset A$. Thus if $A$ is a thick set, it must contain a subset in the form
\begin{align}
\bigcup_{m=1}^\infty\Big\{a_l, a_l+1, a_l+2, \cdots, a_m+m-1\Big\}
\end{align}
for some sequence of integers $a_m\rightarrow\infty$.

It is easy to see that a set $S$ is syndetic$\Leftrightarrow\mathbb{N}\backslash S$ is not thick $\Leftrightarrow S\cap A\ne\emptyset$ for any thick set $A$. Syndetic set and thick set are fundamental concepts in ergodic theory, for details, one may see Furstenberg \cite{fur}.

\begin{definition}
A subset $A\subset\mathbb{N}$ is a Bohr set if there exists a trigonometric polynomial $\psi(t)=\sum_{k=1}^mc_ke^{i\lambda_kt}$, with the $\lambda_k$ real numbers, such that the set
\begin{align}
A'=\{n\in\mathbb{N}:\text{Re}\psi(n)>0\}
\end{align}
is non-empty and $A\supset A'$. When $\psi(0)>0$ we say $A$ is a $Bohr_0$ set.
\end{definition}
As a consequence of the almost periodicity of trigonometric polynomials we can see that a Bohr set is syndetic. We may also define Bohr set and $\text{Bohr}_0$ set in an alternative way, a subset $A\subset\mathbb{N}$ is Bohr set if there exist $m\in\mathbb{N}, \alpha\in\mathbb{T}^m$, and a open set $U\subset\mathbb{T}^m$ such that $\{n\in \mathbb{N}:n\alpha\in U\}$ is contained in $A$; the set $A$ is a $\text{Bohr}_0$ set if additionally $0\in U$.

Bohr-sets are fundamentally abelian in nature. Nowadays it has become apparent that higher order non-abelian Fourier analysis plays an important role both in combinatorial number theory and ergodic theory. Related to this, a higher-order version of Bohr$_0$ sets, namely Nil$_d$ Bohr$_0$-sets, was introduced in \cite{hos}.

\begin{definition}
A subset $A\subset\mathbb{N}$ is a Nil$_d$ Bohr$_0$-set if there exist a $d$-step nilsystem $(X, \mu, T)$, $x_0\in X$, and an open set $U\subset X$ containing $x_0$ such that
\begin{align}
\{n\in\mathbb{N}: T^nx_0\in U\}
\end{align}
is contained in $A$.
\end{definition}

Bergelson, Furstenberg and Weiss introduced the notion of piecewise-Bohr set in \cite{ber}. They defined that a set $A$ is a \emph{piecewise-Bohr} set if $A=S\cap Q$, where $S$ is a Bohr set and $Q$ is a thick set. This notion of piecewise-Bohr set is very simple but weak, a piecewise-Bohr set defined in this manner is even not necessarily syndetic. Then Host and Kra introduced a stronger definition of piecewise-Bohr set, named by \emph{strongly piecewise-Bohr} set in \cite{hos}.
\begin{definition}
The set $A\subset\mathbb{N}$ is said to be strongly piecewise-Bohr, if for every sequence ($J_k: k\ge1$) of intervals whose lengths $|J_k|$ tend to $\infty$, there exists a sequence ($I_j: j\ge1$) of intervals satisfying the following.

\begin{itemize}
  \item[(i)\ \ ] For each $j\ge1$, there exists some $k=k(j)$ such that the interval $I_j$ is contained in $J_k$.
  \item[(ii)~]The lengths $I_j$ tend to infinity.
  \item[(iii)] There exists a Bohr set $B$ such that $B\cap I_j\subset A$ for every $j\ge1$.
\end{itemize}
\end{definition}
Similarly we may define \emph{strongly piecewise-$\text{Bohr}_0$} set. With this definition, both strongly piecewise-Bohr set and strongly piecewise-$\text{Bohr}_0$ set are syndetic.

\section{Zhang-Maynard-Tao's theorem}
\label{sec3}

Let $k$ be a positive integer, we say a given $k$-tuple of integers $H=\{h_1, h_2, \cdots, h_k\}$ is \emph{admissible} if
\begin{align}
\Big|\Big\{n~\text{mod}~p: \prod_{i=1}^k(n+h_i)\equiv0~\text{mod}~p\Big\}\Big|<p, ~\text{for every prime}~p.
\end{align}
In other words, $H$ is admissible if and only if, for any prime $p$, $h_i$'s never occupy all of the residue classes modulo $p$. This is immediately true for all primes $p>k$; so to test this condition for a $k$-tuple of integers $H$ we need only to examine such small primes $p\le k$.

We observe that either Zhang's work or Maynard and Tao's work may follow from a result in the form as
\begin{theorem}
\label{theoremz}
Let $H=\{h_1, h_2, \cdots, h_k\}$ be an $k$-tuple of integers. If $H$ is admissible and $k\ge C$ for some given constant $C>0$, then there are infinitely many integers $n$ such that at least two of the numbers $n+h_1, n+h_2, \cdots, n+h_k$ will be prime.
\end{theorem}
Zhang's Theorem 1 in \cite{zha} proved that $C=3.5\times10^6$ is available in Theorem \ref{theoremz}, then obtained there are infinitely many couples of primes with difference not more than $7\times10^7$. To obtain Maynard-Tao's theorem, they proved that a much smaller value $C=105$ can be used in Theorem \ref{theoremz}. In \cite{pol} they proved $C=50$ is available. If assumed that the primes have level of distribution $\theta$ for every $\theta<1$, Maynard proved that $C$ may take value 5 in the theorem and improved the inferior limit of the difference of primes to 12.

\section{proof of main results}

In this section we prove Theorem \ref{theoremhw} and Corollary \ref{cor}. We begin with an observation that if $A$ is a subset of $\mathbb{N}$ which has more enough elements, then $A$ contains at least an admissible $k$-tuple of integers $H=\{h_1, \cdots, h_k\}$. To find an admissible $k$-tuple of integers, we only need to consider such primes $\mathbb{P}_k=\{p: p\le k\}$. For any prime $p_1\in\mathbb{P}_k$, we have that
\begin{align}
|A|=\sum_{a~mod~p_1}\sum_{{g\in A}\atop{g\equiv a~mod~p_1}}1,
\end{align}
so there exists an integer $b_{p_1}$ such that
\begin{align}
\label{mod}
|\{g\in A: g\equiv b_{p_1}~\text{mod}~p_1\}|\le |A|/p_1.
\end{align}
Let
\begin{align}
A_1=A\setminus \{g\in A: g\equiv b_{p_1}~\text{mod}~p_1\}.
\end{align}
It is easy to see that
\begin{align}
|A_1|\ge|A|\Big(1-\frac1{p_1}\Big)
\end{align}
Then for another prime $p_2\in\mathbb{P}_k$ with $p_2\neq p_1$, we can also choose a $b_{p_2}$ so that
\begin{align}
|\{g\in A_1: g\equiv b_{p_2}~\text{mod}~p_2\}|\le |A_1|/p_2.
\end{align}
The same to $A_1$ we may get a set
\begin{align}
A_2=A_1\setminus \{g\in A_1: g\equiv b_{p_2}~\text{mod}~p_2\}
\end{align}
with
\begin{align}
|A_2|\ge|A_1|\Big(1-\frac1{p_2}\Big)
\end{align}
Repeating this process one prime at a time, with $p$ varying over the elements of $\mathbb{P}_k$, we eventually obtain a set
\begin{align}
\label{2}
A_{\pi(k)}=A\setminus\bigcup_{p\in\mathbb{P}_k}\Big\{g: g\equiv b_p \mod p\Big\}
\end{align}
after $\pi(k)$ steps and the cardinality of this set is
\begin{align}
\label{3}
|A_{\pi(k)}|\ge|A|\prod\limits_{p\in\mathbb{P}_k}\Big(1-\frac1p\Big).
\end{align}
Here $\pi(k)$ denote the number of primes not more than $k$.

From (\ref{2}) we have that, for any $p\le k$, elements of $A_{\pi(k)}$ never occupy all of the residue classes modulo $p$. Thus if
\begin{align}
\label{gek}
|A_{\pi(k)}|\ge k,
\end{align}
any $k$-tuple of integers $H\subset A_{\pi(k)}$ is admissible. By (\ref{3}), to meet the condition (\ref{gek}), we just need assure $A$ large enough that
\begin{align}
A\ge k\prod\limits_{p\le k}\Big(1-\frac1p\Big)^{-1}.
\end{align}
Thus we have that any large enough subset $A$ of $\mathbb{N}$, which has at least $k\prod_{p\le k}\big(1-\frac1p\big)^{-1}$ elements, contains at least an admissible $k$-tuple of integers $H=\{h_1, \cdots, h_k\}$.

However, Theorem \ref{theoremz} tells us that when take the integer $k\ge C$, we may have there are infinitely many integers $n$ such that at least two of the numbers $n+h_1, n+h_2, \cdots, n+h_k$ will be prime. So there must be some integers $h_i,~h_j\in H\subset A$, $h_i>h_j$ that $h_i-h_j$ can be expressed in infinitely many ways as the difference of two primes, that is to say
\begin{align}
\Delta(A)\cap D\neq\emptyset.
\end{align}

From the discussion above we may have the following conclusion. For any subset $A$ of $\mathbb{N}$ with at least $C\prod_{p\le C}\big(1-\frac1p\big)^{-1}$ elements, we have that
\begin{align}
\Delta(A)\cap D\neq\emptyset.
\end{align}
That is to say $D$ is a $\Delta_r^*$-set for any $r$ with
\begin{align}
\label{rge}
     r\ge C\prod_{p\le C}\big(1-\frac1p\big)^{-1}.
\end{align}

Here the constant $C$ is given in Theorem \ref{theoremz}. Unconditionally, the smallest possible value of $C$ that we can take now is 50, to see \cite{pol}. So, from (\ref{rge}), we have $r\ge720.96$, and $D$ is a $\Delta_{721}^*$-set.

If assumed that the primes have level of distribution $\theta$ for every $\theta<1$, Maynard proved in \cite{may} that $C=5$ is available. Thus under this condition, we have $r\ge18.75$ in (\ref{rge}), and then $D$ is a $\Delta_{19}^*$-set.

To prove the corollary, we need the following lemma.
\begin{lemma}
\label{lem}
Every $\Delta^*$-set is a strongly piecewise-Bohr$_0$ set.
\end{lemma}
This lemma is Theorem 2.8 in Host and Kra's work \cite{hos}. We have proved that $D$ is a $\Delta_r^*$-set above, so it is obviously a $\Delta^*$-set. Now it is easy to see that Corollary \ref{cor} is a direct result of Lemma \ref{lem}.

\end{document}